\documentclass[twoside]{amsart}
\usepackage{amssymb,amsmath}
\usepackage[all]{xy}
\usepackage{fancyhdr,color}

\def\bc{\begin{center}}
\def\ec{\end{center}}
\def\no{\noindent}

\def\loa{\longrightarrow}

\def\sul{\sum\limits}

\def\ra{\rightarrow}
\def\ol{\overline}

\def\suse{\subseteq}

\def\sta{\stackrel}
\def\im{{\rm Im}}

\def\fkm{{\frak m}}

\def\ann{{\rm ann}}

\def\bom{\bigotimes}

\def\Ext{{\rm Ext}}

\def\pd{{\rm pd}}

\def\fd{{\rm fd}}

\def\rank{{\rm rank}}

\def\Hom{{\rm Hom}}

\def\lst{\leqslant}
\def\gst{\geqslant}
\def\gpd{G\mbox{-}{\rm pd}}

\def\gl{{\rm gl.dim}}
\def\wgl{w.{\rm gl.dim}}
\def\ggl{G\mbox{-}{\rm gl.dim}}
\def\wggl{w.G\mbox{-}{\rm gl.dim}}
\def\gfd{G\mbox{-}{\rm fd}}

\def\spd{{\rm s.gl.dim}}
\def\gspd{G\mbox{-}{\rm s.gl.dim}}

\def\C{\mathcal C_R}
\def\s{\mathcal S_R}
\def\P{\mathcal P_R}
\def\es{\mathcal S}

\newtheorem{thm}{Theorem}[section]
\newtheorem{defn}[thm]{Definition}
\newtheorem{lem}[thm]{Lemma}
\newtheorem{cor}[thm]{Corollary}
\newtheorem{pro}[thm]{Proposition}
\newtheorem{exa}[thm]{Example}
\newtheorem{rem}[thm]{Remark}

\begin{document}

\vspace*{-2mm}
\bc{\LARGE Super finitely presented modules and Gorenstein projective modules}\ec

\thispagestyle{empty}

\vskip6mm \bc {Fanggui Wang$^1$, Lei Qiao$^2$, and Hwankoo
Kim$^3$}\ec

\vskip2mm 1,2. College of Mathematics, Sichuan Normal University,
Chengdu, Sichuan 610068, China

\bc{wangfg2004@163.com; qiaolei5@yeah.net}\ec

3. School of Computer and Information Engineering, Hoseo University, Asan 336-795, Korea

\bc{hkkim@hoseo.edu}\ec

\vskip6mm
\bc
\begin{minipage}{140truemm}
{\small \no {\bf Abstract}\ \ Let $R$ be a commutative ring. An $R$-module $M$ is said to be super finitely presented if there is an
exact sequence of $R$-modules
$$\cdots \ra P_n\ra\cdots \ra P_1\ra P_0\ra M\ra 0$$
where each $P_i$ is finitely generated projective. In this paper it is shown that if $R$ has the property (B) that every super finitely
presented module has finite Gorenstein projective dimension, then every finitely generated Gorenstein projective module is super
finitely presented. As an application of the notion of super finitely presented modules, we show that if $R$ has the property (C) that every super
finitely presented module has finite projective dimension, then $R$ is $K_0$-regular, i.e., $K_0(R[x_1,\cdots,x_n])\cong K_0(R)$ for all
$n\gst 1$.

\vskip1mm
\no{\bf Keywords:}\ \ Gorenstein projective module, super finitely presented module, super finitely presented dimension,
Gorenstein super finitely presented dimension

\no {\bf MSC(2010):}\ \ 13D02, 13D05, 13D07, 13D15}
\end{minipage}
\ec

\vskip6mm
\section{\bf Introduction}

\vskip3mm
\no Throughout this paper, $R$ denotes a commutative ring with identity and all modules are unitary. For any $R$-module $M$,
$\pd_RM$ and $\fd_RM$ will denote the usual projective and flat dimensions of $M$, respectively. The dual module
$\Hom_{R}(M,R)$ of $M$ is denoted by $M^*$. We use $\gl(R)$ and $\wgl(R)$ to stand for the global dimension and the weak global dimension of $R$,
respectively.

Let $R$ be a Noetherian ring. Following Auslander and Bridger \cite{AB69}, a finitely generated $R$-module $M$ belongs to the
$G$-class $G(R)$, or is called a module of G-dimension zero, if and only if

(1)\ \ $\Ext^{m}_{R}(M,R)=0$ for all $m>0$;

(2)\ \ $\Ext^{m}_{R}(M^*,R)=0$ for all $m>0$; and

(3)\ \ The double duality map $\mu_M:M\ra M^{**}$ is an isomorphism.

This notion had been generalized to an arbitrary rings by Enochs and Jenda in the literature \cite{EJ95}. An $R$-module $M$ is
called \textit{Gorenstein projective} ($G$-projective for short) if there exists an exact sequence of projective $R$-modules
$$\mathbf{P}:=\ \ \cdots\ra P_n\ra\cdots\ra P_1\ra P_0\ra P^0\ra P^1\ra\cdots\ra P^m\ra\cdots \eqno{(1.1)}$$
such that $M\cong \im(P_0\ra P^0)$ and such that $\Hom_R(-,Q)$
leaves the sequence ${\bf P}$ exact whenever $Q$ is a projective
$R$-module. The complex $\mathbf{P}$ is called a {\it complete
projective resolution}. Enochs et al. \cite{EJT93} also introduced
the notion of Gorenstein flat modules. Recall that an $R$-module $M$
is called {\it Gorenstein flat} ($G$-flat for short), if there
exists an exact sequence of flat $R$-modules
$$\mathbf{F}:=\ \ \cdots\ra F_n\cdots\ra F_1\ra F_0\ra F^0\ra F^1\ra\cdots\ra F^m\ra\cdots$$ such
that $M\cong \im(F_0\ra F^0)$ and such that the functor $I\bom_R-$
leaves $\mathbf{F}$ exact whenever $I$ is an injective $R$-module.
Moreover, the Gorenstein projective and flat dimensions of an
$R$-module $M$ are defined in terms of Gorenstein projective and
flat resolutions, respectively, and denoted by $G\mbox{-pd}_RM$ and
$G\mbox{-fd}_RM$, respectively; see \cite{Ho04}. As in the
literature \cite{BM10a}, the Gorenstein global dimension of $R$ is
$$\ggl(R)=\sup\{\gpd_RM\,|\,\mbox{$M$ is an $R$-module}\}.$$

It is well known that every finitely generated projective module is finitely presented. For an integral domain $R$,
Cartier in \cite{Ca58} (also see \cite{Va69}) showed that every projective module of finite rank is finitely generated.
These lead to the following questions:

\vskip1mm
\textbf{Question 1.}\ \ Is every finitely generated $G$-projective module finitely presented?

\vskip1mm
\textbf{Question 2.}\ \ Is every $G$-projective module of finite rank over an integral domain finitely generated?

In order to investigate the two questions above, we use the notion of so-called super finitely presented modules. Recall from \cite{GW14} that an
$R$-module $M$ is said to be {\it super finitely presented} if it admits a projective resolution
$$\cdots \ra P_n\ra\cdots \ra P_1\ra P_0\ra M\ra 0$$
such that each $P_i$ is a finitely generated projective $R$-module.
This notion has received attention in several papers in the
literature. It originated in Grothendieck's notion of a
pseudo-coherent module \cite{BGL71}; in \cite{Bie76, HM09}, the
authors used the term ``$(FP)_\infty$-module'' in the sense of a
super finitely presented module; in \cite{Gl88} it was said to admit
an infinite finite presentation; and in \cite{Be11}, Bennis called
it an infinitely presented module. Using this notion we prove in
Section 3 that if $R$ has the property (B) in the sense that every
super finitely presented module has finite $G$-projective dimension,
then every finitely generated $G$-projective $R$-module $M$ is super
finitely presented and $M^*$ is $G$-projective (Corollary
\ref{c303}); and that if $R$ is a domain having the property (B),
then every $G$-projective $R$-module $M$ of finite rank is also
super finitely presented and $M^*$ is $G$-projective (Corollary
\ref{c307}). Furthermore, to give an example of a ring $R$ that has
the property (B) but is not coherent (see Example \ref{e406}), we
introduce in Section 4 the notions of super finitely presented
dimension and Gorenstein super finitely presented
 dimension of a ring $R$ (Definition \ref{d401}).

Denote $K_0(R)$ the Grothendieck group of a ring $R$ (the definition
will be recalled in Section 5). In algebraic $K$-theory, it is of
interest to determine when a commutative ring $R$ is {\it
$K_0$-regular} (see \cite{La06}), that is, when
$K_0(R[x_1,\cdots,x_n])$ is isomorphic to $K_0(R)$ for all $n$. For
the case that $R$ is Noetherian, the well-known Grothendieck's
Theorem \cite[Chapter II, Theorem 5.8]{La06} (also see
\cite[Corollary 6.2]{Sw68}) states that if every finitely generated
$R$-module has finite projective dimension, then $R$ is
$K_0$-regular. For the coherent case, the results of Quillen
\cite{Qu72,Qu74} suggest that if $R$ is stably coherent (i.e.,
$R[x_1,\dots,x_n]$ is coherent for all $n\geqslant 0$) and if every
finitely presented $R$-module has finite projective dimension, then
$R$ is $K_0$-regular. Now, let us quote from \cite[p. 323]{S70}:
``There are many results of a homological nature which may be
generalized from Noetherian rings to coherent rings. In this
process, finitely generated modules should in general be replaced by
finitely presented modules." More generally, many results of a
homological nature may be generalized from coherent rings to
arbitrary rings if we replace finitely presented modules by super
finitely presented modules (see \cite{GW} for a recent survey). In
line with this point of view, it is also natural to ask the
following question:

\vskip1mm
\textbf{Question 3.}\ \ If every super finitely presented $R$-module has finite projective dimension, then is $R$ $K_0$-regular?

To obtain an affirmative answer to Question 3, we consider the Grothendieck group $S_0(R)$ defined from the family of super finitely
presented $R$-modules. We solve in Section 5 this question by showing that $S_0(R[x])\cong S_0(R)$, and that if $R$ has the property (C) in
the sense that every super finitely presented $R$-module has finite projective dimension, then $S_0(R)\cong K_0(R)$. The method used in this
section is somewhat similar to that of Swan \cite{Sw68}.

Any undefined notions or notation is standard, as in \cite{Ch00,EJ00,Ro79}.

\vskip6mm
\section{\bf On super finitely presented modules}

\vskip3mm \no It is clear that every super finitely presented
$R$-module is finitely presented, but the inverse is not true in
general. In fact, it is easy to show that every finitely presented
$R$-module is super finitely presented if and only if $R$ is a
coherent ring. Therefore, there exists a finitely presented
$R$-module $M$ that is not super finitely presented over any
non-coherent ring $R$. It is well known that every finitely
generated projective $R$-module is super finitely presented.

\begin{lem}\label{l201}
{\rm
Let $0\ra A\ra B\ra C\ra 0$ be an exact sequence of $R$-modules. If any two of $A$, $B$, $C$ are super finitely
presented, then so is the third one.}
\end{lem}

\begin{pro}\label{p202}
{\rm
Let $u\in R$ be neither a zero divisor nor a unit and let $M$ be a super finitely presented $R$-module. If $u$ is regular on $M$, that is, $u$ is a
non-zero-divisor of $M$, then $M/uM$ is a super finitely presented $R/(u)$-module.}
\end{pro}

\no{\bf Proof.}\ \ Let $\cdots \ra F_n\ra\cdots\ra F_1\ra F_0\ra M\ra 0$ be an exact sequence, where all $F_i$ are finitely generated
projective $R$-modules. Since $u$ is regular on $M$, we have the following exact sequence
$$\cdots \ra F_n/uF_n\ra\cdots\ra F_1/uF_1\ra F_0/uF_0\ra M/uM\ra 0.$$
Therefore, $M/uM$ is a super finitely presented $R/(u)$-module.\hfill$\Box$

\begin{pro}\label{p203}
{\rm Let $u\in R$ be neither a zero divisor nor a unit. If $M$ is a super finitely presented $R/(u)$-module, then $M$ is also a super
finitely presented $R$-module.}
\end{pro}

\no{\bf Proof.}\ \ Let $\cdots\ra E_n\ra\cdots\ra E_1\ra E_0\ra M\ra 0$ and $0\ra A\ra Q\ra M\ra 0$ be exact sequences, where all $E_i$
are finitely generated projective $R/(u)$-modules and $Q$ is a finitely generated projective $R$-module. Set
$M_i=\ker(E_i\ra E_{i-1})$ for all $i$ (write $M=M_{-1}$, $Q=Q_{-1}$, and $A=A_{-1}$). Constitute the $3\times 3$ commutative
diagrams for all $i$.
$$\xymatrix@R=16pt@C=28pt{
&0\ar[d]&0\ar[d]&0\ar[d]\\
0\ar[r]&A_i\ar[r]\ar[d]&P_i\ar[d]\ar[r]&A_{i-1}\ar[r]\ar[d]&0\\
0\ar[r]&Q_i\ar[d]\ar[r]&F_i\ar[r]\ar[d]&Q_{i-1}\ar[d]\ar[r]&0\\
0\ar[r]&M_i\ar[r]\ar[d]&E_i\ar[r]\ar[d]&M_{i-1}\ar[r]\ar[d]&0\\
&0&0&0}$$
where all $Q_i$ are finitely generated projective $R$-modules. Because each $M_i$ is a finitely presented
$R/(u)$-module, $M_i$ is a finitely presented $R$-module by \cite[Theorem 5.1.5]{Wa06}. Since $\pd_RE_i=1$, $P_i$ is finitely
generated projective for all $i$. Hence we have the following exact sequence
$$\cdots\ra P_n\ra\cdots\ra P_1\ra P_0\ra A\ra 0.$$
Consequently, $A$ is a super finitely presented $R$-module. Therefore, $M$ is also a super finitely presented $R$-module by
Lemma \ref{l201}.\hfill$\Box$

\vskip2mm Let $M$ be an $R$-module. Then $M\bigotimes_R R[x]$ is an
$R[x]$-module which we write as $M[x]$. Thus, elements of $M[x]$ are
of the form $\sul^n_{i=0}m_i\otimes x^i$, where $i\geqslant 0$ and
$m_i\in M$. The proofs of the following two propositions are easy,
so we omit the proofs of them.

\begin{pro}\label{p204}
{\rm Let $R$ and $T$ be rings, and $\varphi:R\ra T$ be a
homomorphism. Suppose $T$ is a flat $R$-module. If $M$ is a super
finitely presented $R$-module, then $T\bom_RM$ is a super finitely
presented $T$-module. In particular, if $M$ is a super finitely
presented $R$-module, then the polynomial module $M[x]$ is a super
finitely presented $R[x]$-module.}
\end{pro}

\begin{pro}\label{p205}
{\rm
Let $R$ and $T$ be rings, $\varphi:R\ra T$ be a homomorphism, and $T$ be a finitely generated projective $R$-module. If $M$ is a
super finitely presented $T$-module, then $M$, as an $R$-module, is super finitely presented.}
\end{pro}

\begin{pro}\label{p206}
{\rm
Let $M$ be a finitely generated $G$-projective $R$-module. Then we have:

(1)\ \ There is an exact sequence of $R$-modules
$$0\ra M\ra P_0\ra P_1\ra \cdots \ra P_m\ra \cdots,\eqno{(2.1)}$$
in which each $P_i$ is finitely generated projective, and each cosyzygy of (2.1) is a finitely generated $G$-projective;

(2)\ \ $M^*$ is super finitely presented.}
\end{pro}

\no{\bf Proof.}\ \ See \cite[Lemma 2.3 and Corollary 2.4]{GW14}.\hfill$\Box$

\vskip6mm
\section{\bf Finitely generated Gorenstein projective modules}

\vskip3mm \no Recall that an $R$-module $X$ is called torsionless if
the natural map $\mu:X\ra X^{**}$ is monomorphic and that $X$ is
called reflexive if $\mu$ is isomorphic. It is well-known that $X$
is torsionless if and only if $X$ can be imbedded in a product of
copies of $R$. Therefore, all submodules of a free module are
torsionless.

\begin{lem}\label{l301}
{\rm Let $0\ra M\ra P_0\ra C\ra 0$ be an exact sequence, where $P_0$ is finitely generated projective, and $M$ and $C$ are
$G$-projective. Then we have:

(1)\ \ $M$ and $M^*$ are reflexive.

(2)\ \ $\Ext_R^i(M^*,R)=0$ for all $i>0$.}
\end{lem}

\no{\bf Proof.}\ \ (1)\ \ Since $C$ is $G$-projective, we have the exact sequence $0\ra C^*\ra P_0^*\ra M^*\ra 0$. Taking duality again we
obtain the following commutative diagram with exact rows:
$$\xymatrix@R=17pt@C=30pt{
0\ar[r]&M\ar[r]\ar[d]_{\mu_M}&P_0\ar[r]\ar[d]^{\cong}&C\ar[r]\ar[d]^{\mu_C}&0\\
0\ar[r]&M^{**}\ar[r]&P_0^{**}\ar[r]&C^{**}\ar[r]&\Ext_R^1(M^*,R)\ar[r]&0}$$
where $\mu_X:X\ra X^{**}$ is the natural homomorphism for $X=M,C$. Since $M$ and $C$ are torsionless, we have that $\mu_M$ is
isomorphic by Five Lemma. Thus $M$ is reflexive. Taking the duality we have that $M^*$ is also reflexive.

(2)\ \ From the commutative diagram above we have also that $\mu_C$ is isomorphic by the same argument. Consequently,
$\Ext_R^1(M^*,R)=0$. Since $C$ is finitely generated $G$-projective, we have by applying Proposition \ref{p206}(1) on $C$ the following
exact sequence
$$0\ra M\ra P_0\ra P_1\ra \cdots\ra P_n\ra\cdots ,\eqno{(3.1)}$$
where all $P_i$ are finitely generated projective. Hence we obtain the exact sequence
$$\cdots\ra P_n^*\ra \cdots \ra P_1^*\ra P_0^*\ra M^*\ra 0.\eqno{(3.2)}$$
Let $C_i$ be the $i$-th cosyzygy ($C_0=C$) of (3.1). Then $C_i^*$ is the $i$-th syzygy of (3.2). Thus we have by applying the argument above
that $C_i$ is reflexive and $\Ext_R^1(C_i^*,R)=0$. Consequently, we have $\Ext_R^i(M^*,R)=0$ for all $i>0$.\hfill$\Box$

\vskip2mm For convenience we say that all kernels of arrows in the
exact sequence (1.1) are syzygies.

\begin{thm}\label{t302}
{\rm
Let $M$ be a finitely generated $G$-projective $R$-module. If $\gpd_RM^*<\infty$, then $M$ is super finitely presented and $M^*$ is
$G$-projective.}
\end{thm}

\vskip1mm \no{\bf Proof.}\ \ Let $0\ra A\ra P\ra M\ra 0$ be an exact sequence, where $P$ is a finitely generated projective $R$-module
and $A$ is a $G$-projective $R$-module. Then $0\ra M^*\ra P^*\ra A^*\ra 0$ is an exact sequence and $A^*$ is reflexive by Lemma
\ref{l301}. Now we let $\gpd_RM^*\lst n+1<\infty$ by hypothesis. Then, by Proposition \ref{p206}(2), there is an exact sequence
$$0\ra X\ra F_n\ra \cdots \ra F_1\ra F_0\ra P^*\ra A^*\ra 0$$
in which all $F_i$ are finitely generated projective and $X$ is finitely generated $G$-projective. By Proposition \ref{p206}(1),
there exists an exact sequence
$$0\ra X\ra Q_n\ra \cdots \ra Q_1\ra Q_0\ra Q\ra N\ra 0,$$
such that $Q$ and all $Q_i$ are finitely generated projective modules and $N$ is finitely generated $G$-projective. Because all
$\im (Q_i\ra Q_{i-1})$ and all $\im(Q_0\ra Q)$ are $G$-projective, we have the following commutative diagram with exact rows:
$$\xymatrix@R=15pt@C=25pt{
0\ar[r]&X\ar@{=}[d]\ar[r]&Q_n\ar[d]\ar[r]&{\cdots}\ar[r]&Q_1\ar[d]\ar[r]&Q_0\ar[d]\ar[r]&Q\ar[r]\ar[d]&N\ar[d]\ar[r]&0\\
0\ar[r]&X\ar[r]&F_n\ar[r]&{\cdots}\ar[r]&F_1\ar[r]&F_0\ar[r]&P^*\ar[r]&A^*\ar[r]&0}$$
Thus we obtain the following exact sequence
$$\mbox{$0\ra Q_n\ra Q_{n-1}\bigoplus F_n\ra \cdots \ra Q_0\bigoplus F_1\ra Q\bigoplus F_0\ra N\bigoplus P^*\ra A^*\ra 0$}.\eqno{(3.3)}$$
Let $Y_i$ be the $i$-th syzygy ($Y_{n+1}=Q_n$) of exact sequence (3.3). By Lemma \ref{l301},
$\Ext_R^1(Y_i,R)\cong\Ext_R^{i+1}(A^*,R)=0$. Since $A^*$ and all $Y_i$ are super finitely presented, we have $\Ext_R^1(A^*,F)=0$ and
$\Ext_R^1(Y_i,F)=0$ for any projective module $F$. Since $0\ra Q_n\ra Q_{n-1}\bigoplus F_n\ra Y_n\ra 0$ is exact, and
$\Ext_R^1(Y_n,Q_n)=0$, we have that this sequence is splitting, whence $Y_n$ is projective. By repeating this process we obtain that
$Y_{n-1},\cdots,Y_1,Y_0$ are projective. Now we consider the exact sequence $0\ra Y_0\ra N\bigoplus P^*\ra A^*\ra 0$. We have also that
this sequence is splitting, whence we have
$$\mbox{$(Y_0)^*\bigoplus A=(Y_0)^*\bigoplus A^{**}\cong N^*\bigoplus P$}.$$
Because $N$ is a finitely generated $G$-projective, $N^*$ is finitely generated by Proposition \ref{p206}. Thus $A$ is finitely
generated. Therefore, $M$ is finitely presented.

From the exact sequence $0\ra M^*\ra P^*\ra A^*\ra 0$, we have $\gpd_RA^*\lst \gpd_RM^*+1<\infty$. By the same argument we have that
$A$ is finitely presented. Continuing this process we can obtain that $M$ is super finitely presented.

Since $M$ is super finitely presented, $M$ has a complete projective resolution (1.1) in
which all $P_i$ and $P^j$ are finitely generated projective. Thus we have the exact sequence
$$\cdots \ra (P^m)^*\ra \cdots (P^1)^*\ra (P^0)^*\ra (P_0)^*\ra (P_1)^*\ra\cdots\ra (P_n)^*\ra\cdots\eqno{(3.4)}$$
in which $M^*\cong \im((P^0)^*\ra (P_0)^*)$. Let $X$ be any syzygy of (1.1). Then $X$ is super finitely generated $G$-projective and
$X^*$ is a syzygy of (3.4). By Lemma \ref{l301}, $\Ext_R^1(X^*,R)=0$, and whence $\Ext_R^1(X^*,Q)=0$ for any
projective module $Q$ because $X^*$ is super finitely presented by Proposition \ref{p206}. Therefore, $M^*$ is
$G$-projective.\hfill$\Box$

\vskip2mm To determine rings for which every finitely generated
$G$-projective module is finitely presented, we say that a ring $R$
has the property (B) if every super finitely presented $R$-module
has finite $G$-projective dimension. Certainly, if $\ggl(R)<\infty$,
then $R$ has the property (B). By Theorem \ref{t302}, we have the
following corollaries.

\begin{cor}\label{c303}
{\rm
Let $R$ have the property (B) and let $M$ be a finitely generated $G$-projective $R$-module. Then $M$ is super finitely
presented and $M^*$ is $G$-projective.}
\end{cor}

\begin{cor}\label{c304}
{\rm
If $\ggl(R)<\infty$, then every finitely generated $G$-projective $R$-module $M$ is super finitely presented and $M^*$ is
$G$-projective.}
\end{cor}

Now we investigate Question 2. Let $R$ be an integral domain with quotient field $K$. The rank of an $R$-module $M$ is defined by
$\rank(M)=\dim_K(K\bom_RM)$. It is well-known that $M$ is a torsion module if and only if $\rank(M)=0$.

\begin{lem}\label{l305}
{\rm
Let $R$ be an integral domain with quotient field $K$ and $M$ be a $G$-projective $R$-module of finite rank. Then there
exists an exact sequence $0\ra M\ra P\ra C\ra 0$ such that $P$ is finitely generated free and $C$ is $G$-projective.}
\end{lem}

\no{\bf Proof.}\ \ Since $M$ is $G$-projective, there is an exact sequence $0\ra M\sta{\varphi}{\ra} F\ra N\ra 0$ in which $F$ is
projective and $N$ is $G$-projective. Without loss of generality we assume that $F$ is free with a basis $\{e_i\,|\,i\in \Gamma\}$.
Clearly, $M$ is torsion-free. Let $\rank(M)=n$. Suppose $x_1,\cdots,x_n\in M$ is a basis of $K\bom_RM$ over $K$. Thus there
exists a finite subset $\{e_1,\cdots,e_m\}$ such that
$$\varphi(x_k)=\sul_{j=1}^ma_{kj}e_j,\qquad\qquad a_{kj}\in R,\qquad k=1,2,\cdots,n.$$
Set $P=\bigoplus\limits_{j=1}^mRe_j$ and $F_1=\bigoplus\limits_{i\not=1,\cdots,m}Re_i$. Then $F=P\bigoplus F_1$, $P$ is finitely
generated free, and $F_1$ is free.

If $x\in M$, then we have $\varphi(x)=\sul_{i\in \Gamma}r_ie_i$, $r_i\in R$, where $\varphi$ has finite support. Set
$x=\frac{b_1}{s}x_1+\cdots+\frac{b_n}{s}x_n$, for $b_k,s\in R$ with $s\not=0$. Hence we get
$$\varphi(sx)=\sul_{i\in \Gamma}sr_ie_i=\sul_{k=1}^n\sul_{j=1}^ma_{kj}b_ke_j.$$
Consequently, if $i\in\Gamma$ with $i\not=1,\cdots,m$, then $r_i=0$.
Hence $\varphi(M)\suse P$. Set $C=P/\varphi(M)$. Thus $0\ra M\ra P\ra C\ra 0$ is exact and we can get the following commutative
diagram with exact rows
$$\xymatrix@R=15pt@C=25pt{
0\ar[r]&M\ar@{=}[d]\ar[r]&P\ar[d]\ar[r]&C\ar[d]\ar[r]&0\\
0\ar[r]&M\ar[r]&F\ar[r]&N\ar[r]&0}$$ Consequently, $0\ra C\ra N\ra
F_1\ra 0$ is exact by Snake Lemma, and hence have $N\cong
F_1\bigoplus C$. Therefore, $C$ is $G$-projective.\hfill$\Box$

\begin{thm}\label{t306}
{\rm
Let $R$ be an integral domain and let $M$ be a $G$-projective $R$-module of finite rank. If $\gpd_RM^*<\infty$, then $M$ is super finitely
presented and $M^*$ is $G$-projective.}
\end{thm}

\no{\bf Proof.}\ \ By Lemma \ref{l305} we have an exact sequence $0\ra M\ra P\ra C\ra 0$, where $P$ is finitely generated projective
and $C$ is finitely generated $G$-projective. Since $\gpd_RM^*<\infty$, we have $\gpd_RC^*<\infty$. By Theorem
\ref{t302}, $C$ is super finitely presented. Hence $M$ is super finitely generated. By using Theorem \ref{t302} again,
$M^*$ is $G$-projective.\hfill$\Box$

\begin{cor}\label{c307}
{\rm
If $R$ is an integral domain having the property (B), then every $G$-projective $R$-module $M$ of finite rank is super finitely
presented and $M^*$ is $G$-projective.}
\end{cor}

\begin{cor}\label{c308}
{\rm
If $R$ is an integral domain with $\ggl(R)<\infty$, then every $G$-projective $R$-module $M$ of finite rank is super finitely presented
and $M^*$ is $G$-projective.}
\end{cor}

\vskip2mm In \cite{Be11} Bennis proved that a super finitely
generated module $M$ is $G$-projective if and only if $M$ is
$G$-flat. Thus we have the following corollaries.

\begin{cor}\label{c309}
{\rm
If $R$ has the property (B), then every finitely generated $G$-projective $R$-module is $G$-flat.}
\end{cor}

\begin{cor}\label{c310}
{\rm
If $\ggl(R)<\infty$, then every finitely generated $G$-projective module is $G$-flat.}
\end{cor}

In \cite{BM09} it is proved that if $M$ is an $R$-module, then $\gpd_{R[x]}M[x]\lst \gpd_RM$; and if all projective
$R[X]$-modules have finite injective dimension, then $\gpd_{R[x]}M[x]=\gpd_RM$. In fact, this equality is true under any
case.

\begin{lem}\label{l311}
{\rm
Let $M$ be an $R$-module. Then $\gpd_{R[x]}{M[x]}=\gpd_RM$. Moreover, an $R$-module $M$ is $G$-projective if and only if $M[x]$
is a $G$-projective $R[x]$-module.}
\end{lem}

\no{\bf Proof.}\ \ By \cite[Lemma 2.5]{BM09}, we have $\gpd_{R[x]}{M[x]}\lst\gpd_RM$. So the equality holds for the case
$\gpd_{R[x]}{M[x]}=\infty$.

Now suppose that $\gpd_{R[x]}{M[x]}\lst n<\infty$. Let $0\ra F_n\ra \cdots\ra F_1\ra F_0\ra M[x]\ra 0$ be an exact sequence, where
$F_0,\ldots,F_{n-1}$ are free $R[x]$-modules. Thus $F_n$ is a $G$-projective $R[x]$-module. Note that $x$ is regular on
$F_0,\ldots,F_n$ and on $M[x]$. Hence we get the following exact sequence
$$0\ra F_n/xF_n\ra \cdots\ra F_1/xF_1\ra F_0/xF_0\ra M\ra 0. \eqno{(3.5)}$$
By \cite[Theorem 3.1]{BM10a}, every $F_i/xF_i$ is a $G$-projective $R$-module. Hence the sequence (3.5) is a $G$-projective resolution
of $M$ as an $R$-module. Thus $\gpd_RM\lst n$. Therefore, $\gpd_{R[x]}{M[x]}=\gpd_RM$.\hfill$\Box$

\vskip2mm The following lemma is \cite[Lemma 1]{WC09}. But the
article reference is in Chinese, so we include a proof for the
convenience of the reader.

\begin{lem}\label{lemma}{\rm Given a commutative diagram in
the category of all $R$-modules with exact rows, $$\xymatrix{
  0\ar[r]&A\ar@{^{(}->}[r]\ar@{^{(}->}[d]&B\ar[r]^f\ar@{^{(}->}[d]&C\ar[r]\ar[d]^h&0\\
0\ar[r]&A_1\ar@{^{(}->}[r]&B_1\ar[r]^g&C_1\ar[r]&0, }$$ where each
``$\hookrightarrow$'' is an inclusion. Then $A=A_1\bigcap B$ if and
only if $h$ is a monomorphism.}
\end{lem}

\no{\bf Proof.}\ \ Suppose that $h$ is a monomorphism. For each
$y\in A_1\bigcap B$, we have $g(y)=hf(y)=0$, and so $f(y)=0$, i.e.,
$y\in A$. It follows that $A=A_1\bigcap B$. Conversely, let
$A=A_1\bigcap B$, and let $h(z)=0$ for some $z\in C$. Since $f$ is
epimorphic, there is a $y\in B$ with $f(y)=z$. Thus,
$g(y)=hf(x)=h(z)=0$, and so $y\in A_1\bigcap B=A$. Hence,
$z=f(y)=0$. This proves that $h$ is a monomorphism.\hfill$\Box$

\vskip2mm
Now we show that the property (B) is stable under polynomial ring extensions.

\begin{lem}\label{l312}
{\rm
Let $0\ra M\ra F[x]\ra N\ra 0$ be an exact sequence over $R[x]$, where $F$ is a finitely generated free $R$-module and $N$ is a super
finitely presented submodule of a free $R[x]$-module. Then there is an exact $R[x]$-sequence
$$0\loa A[x]\loa B[x]\loa M\loa 0\eqno{(3.6)}$$
in which $A$ and $B$ are super finitely presented $R$-modules.}
\end{lem}

\no{\bf Proof.}\ \ Note that $M$ is a super finitely presented $R[x]$-module. Let $\{z_1,\cdots,z_s\}$ be a generating set of $M$ and set
$F_k=F+Fx+\cdots+Fx^{k-1}$ for all $k\gst 1$. Then each $F_k$ are finitely generated free $R$-module. Write $B=M\bigcap F_k$. Then we
have the following commutative diagram with exact $R$-rows:
$$\xymatrix@R=16pt@C=27pt{
0\ar[r]&B\ar[r]\ar[d]&F_k\ar[r]^{\pi}\ar[d]&F_k/B\ar[r]\ar[d]^{\alpha}&0\\
0\ar[r]&M\ar[r]&F[x]\ar[r]&N\ar[r]&0}$$ where the homomorphism
$\alpha$ comes from diagram chasing. By Lemma \ref{lemma}, $\alpha$
is a monomorphism. Set $T=R[x]/(x^k)$. Because $N$ is a submodule of
a free $R[x]$-module, $x^k$ is regular on $N$, whence we have the
following exact $T$-sequence
$$0\loa M/x^kM\loa F[x]/x^kF[x]\loa N/x^kN\loa 0.$$
Consequently, we have the following commutative diagram with exact rows:
$$\xymatrix@R=16pt@C=27pt{
0\ar[r]&B\ar[r]\ar[d]_{\beta}&F_k\ar[r]\ar[d]^{\cong}&F_k/B\ar[r]\ar[d]^{\gamma}&0\\
0\ar[r]&M/x^kM\ar[r]&F[x]/x^kF[x]\ar[r]&N/x^kN\ar[r]&0}$$ Note that
the middle vertical arrow is $R$-isomorphic. Thus the composition homomorphism $\beta:B\ra M\ra M/x^kM$ is an $R$-monomorphism.

Take a sufficient large $k$ such that $z_1,\cdots,z_s\in B$. For such a $k$ we have that $\beta$ is an $R$-epimorphism, and hence
$\beta$ is an $R$-isomorphism. Therefore, $\gamma$ is also an $R$-isomorphism. By Proposition \ref{p202}, $M/x^kM$ is a super
finitely presented $T$-module. By Proposition \ref{p205}, $M/x^kM$ is a super finitely presented $R$-module. Therefore, $B$ is a super
finitely presented $R$-module.

Set $A=M\bigcap F_{k-1}$. Then $A\suse B$. Since $A=\{u\in M\,|\,xu\in B\}$, we can define $\psi:B[x]\ra M$ by
$$\psi(\sul_{i=0}^nb_i\otimes x^i)=\sul_{i=0}^nx^ib_i,\qquad b_i\in B.$$
Then $\psi$ is an $R[x]$-homomorphism. Since $z_1,\cdots,z_s\in B$, we have $\psi$ is an epimorphism.

Define $\varphi:A[x]\ra B[x]$ by
$$\varphi(\sul_{i=0}^na_i\otimes x^i)=\sul_{i=0}^na_i\otimes x^{i+1}-\sul_{i=0}^n(xa_i)\otimes x^i,\qquad a_i\in A.$$
Then $\varphi$ is also an $R[x]$-homomorphism. It is routine to show
that $0\loa A[x]\sta{\varphi}{\loa} B[x]\sta{\psi}{\loa} M\ra 0$ is
an exact sequence over $R[x]$. By Proposition \ref{p204}, $B[x]$ is
a super finitely presented $R[x]$-module. By Lemma \ref{l201},
$A[x]$ is a super finitely presented $R[x]$-module. Thus $A$ is a
super finitely presented $R$-module by Proposition
\ref{p202}.\hfill$\Box$

\begin{thm}\label{t313}
{\rm
If $R$ has the property (B), then so does the polynomial ring $R[x]$.}
\end{thm}

\no{\bf Proof.}\ \ Let $N$ be a super finitely presented $R[x]$-module. Then there is an exact $R[x]$-sequence $0\ra M\ra
F[x]\sta{g}{\ra} N\ra 0$, where $F$ is a finitely generated free $R$-module and $M$ is a super finitely presented $R[x]$-module. Note
that $\gpd_{R[x]}N<\infty$ if and only if $\gpd_{R[x]}M<\infty$. Replacing $N$ with $M$ we can assume without loss of generality that $N$ is a
submodule of some finitely generated free $R[x]$-module. By Lemma \ref{l312} we have an exact sequence $0\ra A[x]\ra B[x]\ra M\ra 0$, where
$A$ and $B$ are super finitely presented $R$-modules. Since $\gpd_RA<\infty$ and $\gpd_RB<\infty$ by hypothesis, we have that
$\gpd_{R[x]}A[x]<\infty$ and $\gpd_{R[x]}B[x]<\infty$ by Lemma
\ref{l311}. By \cite[Lemma 2.4]{BM10}, $\gpd_{R[x]}M<\infty$, and hence $\gpd_{R[x]}N<\infty$. Therefore, $R[x]$ has the property
(B).\hfill$\Box$

\vskip6mm
\section{\bf The super homological dimensions of a ring}

\vskip3mm \no To exhibit an example of a ring having the property
(B) that is not coherent, we are in the position to define the the
super finitely presented dimension and Gorenstein super finitely
presented dimension of a ring $R$.

\vskip2mm
\begin{defn}\label{d401}
{\rm
For a ring $R$, its super finitely presented dimension is defined by
$$\spd(R)=\sup\{\pd_RM\,|\,\mbox{$M$ is a super finitely presented $R$-module}\},$$
and its Gorenstein super finitely presented dimension is defined by
$$\gspd(R)=\sup\{\gpd_RM\,|\,\mbox{$M$ is a super finitely presented $R$-module}\}.$$}
\end{defn}

Clearly, we have $\gspd(R)\lst \spd(R)$, and if $\spd(R)<\infty$, then $\gspd(R)=\spd(R)$ by \cite[Proposition 2.27]{Ho04}. Moreover,
$\spd(R)\lst\gl(R)$ and $\gspd(R)\lst \ggl(R)$. If $R$ is a coherent ring, then $\spd(R)=\wgl(R)$.

\vskip2mm
Certainly, if $\gspd(R)<\infty$, then $R$ has the property (B).

\vskip2mm
In \cite{BM10a}, the weak Gorenstein global dimension of a ring $R$, denoted by $\wggl(R)$, is defined as
$\wggl(R)=\mbox{sup}\{\mbox{\gfd}_R(M)|\mbox{$M$ is an $R$-module}\}$.

\begin{pro}\label{p402}
{\rm
For any ring $R$, we have:

(1)\ \ $\spd(R)\lst \wgl(R)$.

(2)\ \ $\gspd(R)\lst \wggl(R)$.}
\end{pro}

\no{\bf Proof.}\ \ (1)\ \ It is easy from the fact that finitely presented flat modules are projective.

(2)\ \ It follows directly from the result of Bennis \cite{Be11}
that super finitely presented $G$-flat modules are
$G$-projective.\hfill$\Box$

\begin{pro}\label{p403}
{\rm
Let $u\in R$ be a non-zero-divisor and nonunit element.

(1)\ \ If $\spd(R/(u))<\infty$, then $\spd(R/(u))+1\lst \spd(R)$.

(2)\ \ $\gspd(R/(u))+1\lst \gspd(R)$. Therefore, If $\gspd(R)=n<\infty$, then $\gspd(R/(u))\lst n-1$.}
\end{pro}

\no{\bf Proof.}\ \ (1)\ \ Let $\spd(R/(u))=n$. Then there is a super
finitely presented $R/(u)$-module $M$ with $\pd_{R/(u)}M=n$. Thus
$\pd_RM=n+1$ by \cite[Theorem 3, Part III]{K69}. It implies
$\spd(R/(u))+1\lst \spd(R)$.

(2)\ \ Let $M$ be a super finitely presented $R/(u)$-module. Then $M$ is a super finitely presented $R$-module by Proposition
\ref{p203}. By \cite[Theorem 4.1]{BM10}, $\gpd_{R/(u)}M+1=\gpd_RM\lst\gspd(R)$.\hfill$\Box$

\begin{lem}\label{l404}
{\rm
Let $R$ be a domain. If $R/(u)$ is coherent for any nonzero and nonunit $u\in R$, then $R$ is coherent.}
\end{lem}

\no{\bf Proof.}\ \ Let $I$ be a finitely generated proper ideal of $R$. Take $u\in I$ with $u\not=0$. Since
$R/(u)$ is coherent, $I/(u)$ is finitely presented over $R/(u)$. Hence $I/(u)$ is finitely presented over $R$ by \cite[Theorem 5.1.5]{Wa06}
(or by Proposition \ref{p203}). Since $0\ra (u)\ra I\ra I/(u)\ra 0$ is exact, $I$ is finitely presented. Consequently, $R$ is coherent.\hfill$\Box$

\begin{lem}\label{l405}
{\rm
Let $R$ be a ring and let $a\in R$ such that $a\not=0$ and $a^2=0$. Set $I=(a)$. If $\ann(I)=I$, then $I$ is super finitely presented
and $\pd_RI=\infty$. Therefore, $\spd(R)=\infty$, and hence $\wgl(R)=\infty$. }
\end{lem}

\no{\bf Proof.}\ \ Because $\ann(I)=I$, the sequence $0\ra I\ra R\ra I\ra 0$ is exact. Thus the sequence
$$\cdots\ra R\ra R\ra \cdots \ra R\ra R\ra I\ra 0$$
is exact and all syzygies of this sequence are $I$. Therefore, $I$ is super finitely presented, and $\pd_RI<\infty$ if and only if $I$ is projective.
Hence we are done by showing that $I$ is not projective.

Since $I$ is projective if and only if $I_\fkm$ is free for every maximal ideal $\fkm$ of $R$, we assume without loss
of generality that $R$ is local. If $I$ is projective, then $I$ is free and $a$ must be a basis of $I$, but this will
contradict the fact $a^2=0$.\hfill$\Box$

\vskip2mm
\begin{exa}\label{e406}
{\rm
Now we exhibit a ring having the property (B) that is not coherent and $\wgl(R)=\infty$. In \cite{So68} Soublin gave a ring $S$ that is coherent
with $\wgl(S)=2$, but $S[x]$ is not coherent. By Alfonsi's Reduction Theorem(see \cite[Theorem 7.2.6]{Gl89}), there is a maximal ideal $\fkm$ of $S$ such that
$S_\fkm[x]$ is not coherent. Thus $S_\fkm$ is a GCD domain by \cite[Corollary 6.2.10]{Gl89}. Set $D=S_\fkm[x]$. Then $D$ is not coherent with
$\wgl(D)=3$. By Lemma \ref{l404}, there is a nonzero and nonunit $u\in D$ such that $D/(u)$ is not coherent. Set $R=D/(u^2)$.
By \cite[Theorem 4.1.1(1)]{Gl89}, $R$ is not coherent. By Proposition \ref{p403}, $\gspd(R)\lst 2$. Hence $R$ has the property (B).

Write $a=\ol{u}$ and $I=(a)$. We claim that $\ann(I)=I$. Obviously,
$a\not=0$ and $a^2=0$. Hence, $I\subseteq\ann(I)$. On the other
hand, let $r\in R$ with $ra=0$ and denote $r=\ol{d}$ for some $d\in
D$. Then $du\in (u^2)$. Thus, $d\in (u)$, and so $r\in (a)=I$.
Therefore, $\ann(I)=I$. By Lemma \ref{l405}, $\wgl(R)=\infty$.}
\end{exa}

\vskip6mm
\section{\bf On $K_0$-regularity of rings}

\vskip3mm
\no In this section, we give an affirmative answer to Question 3 mentioned in Introduction. We follow the clue of the so-called
``Grothendieck construction'' which can be found in \cite{LS75}.

For an $R$-module $M$ in a given family $\C$ of $R$-modules, let $(M)$ denote the isomorphism class of $M$. Let $G$ be the
free abelian group on the basis $\{(M):M\in \C\}$, and let $H$ be the subgroup generated by all elements of $G$
of the form $(M)-(M_1)-(M_2)$ whenever
$$0\ra M_1\ra M\ra M_2\ra 0\eqno{(5.1)}$$
is an exact sequence in $\C$. Then the {\it Grothendieck group} of $\C$, denoted by $K_0(\C)$, is defined as the
quotient group $G/H$. For $M\in\C$, the image of $(M)$ in $K_0(\C)$ will be denoted by $[M]$ (or, if necessary, $[M]_{\C}$). Hence,
whenever we have an exact sequence (5.1) in $\C$, the relation $[M]=[M_1]+[M_2]$ holds in $K_0(\C)$. Moreover, it is easy
to see that the group $K_0(R)$ satisfies the so-called ``universal property'', which can be described as follows. Let $(L,+)$ be any
abelian group and let $f:\C\ra L$ be a map such that
\begin{enumerate}
\item for $M\in\C$, the image $f(M)$ depends only on the isomorphism class of $M$;
\item for each exact sequence (5.1) in $\C$, we have $f(M)=f(M_1)+f(M_2)$.
\end{enumerate}
Then there exists a unique group homomorphism $h:K_0(\C)\ra L$ such that $h([M])=f(M)$ for all $M\in\C$.

\begin{rem}\label{r501}
{\rm (1)\ \ Let $\P$ be the family of all finitely generated
projective $R$-modules, it is well-known that the Grothendieck group
of $\P$ is denoted simply by $K_0(R)$, which is called the
Grothendieck group of $R$. As we mentioned in Introduction, a
commutative ring $R$ is called $K_0$-regular if
$K_0(R[x_1,\cdots,x_n])$ is isomorphic to $K_0(R)$ for all $n$.

(2)\ \ If $\C$ is the family of all finitely generated $R$-modules, then the Grothendieck group of $\C$ is
denoted by $G_0(R)$, see \cite{Sw68}.

(3)\ \ If $\C$ is the family of all finitely presented $R$-modules, then the Grothendieck group of $\C$ is
denoted by $K_0({\rm Mod}fp(R))$, see \cite{Ge74}.

(4)\ \ For our purpose, we use $\s$ to denote the family of all super finitely presented $R$-modules and use $S_0(R)$ to
denote the Grothendieck group of $\s$.}
\end{rem}

For convenience, we say that a ring $R$ has the property (C) if every super finitely presented $R$-module has finite projective dimension.
Obviously, if $\spd(R)<\infty$, then $R$ has the property (C). Also, by Proposition \ref{p402}, every ring with finite weak global dimension
has the property (C).

\vskip2mm
Note that $\P\suse\s$. Hence it is obvious that the inclusion of $\P$ into $\s$ induces a group
homomorphism $\theta_R:K_0(R)\ra S_0(R)$ by defining $\theta_R([M])=[M]$ for all $M\in\P$. Note that $\theta_R$ is not
a monomorphism in general. But we have the following:

\begin{pro}\label{p501}
{\rm
If $R$ has the property (C), then the natural group homomorphism $\theta_R:K_0(R)\ra S_0(R)$ is an isomorphism.}
\end{pro}

\no{\bf Proof.}\ \ Let $M\in\s$. Then $M$ admits a finite projective resolution $0\ra P_n\ra\cdots\ra P_1\ra P_0\ra M\ra 0$,
where $P_i\in\P$. Define $f:\s\ra K_0(R)$ by $f(M)=\sum\limits^n_{i=0}(-1)^i[P_i]$. Then by the generalized Schanuel Lemma
(cf.\cite[Exercise 3.37]{Ro79}), we see that $f$ is a well-defined map. And if $0\ra M_1\ra M\ra M_2\ra 0$ is an exact sequence
in $\s$, then, by the Horseshoe Lemma (cf.\cite[Lemma 6.20]{Ro79}), we have $f(M)=f(M_1)+f(M_2)$. Therefore, $f$ induces
a homomorphism $\nu:S_0(R)\ra K_0(R)$ such that $\nu([M])=f(M)=\sum\limits^n_{i=0}(-1)^i[P_i]$ for $M\in\s$.
It is routine to check that $\nu\theta_R={\bf 1}_{K_0(R)}$ and $\theta_R\nu={\bf 1}_{S_0(R)}$, which implies that $\theta_R$ is an
isomorphism.\hfill$\Box$

\vskip2mm Let $\varphi:R\ra T$ be a ring homomorphism such that $T$
as an $R$-module is flat. Define $f:\s\ra S_0(T)$ by $f(M)=[T\bom_R
M]$ for $M\in \s$. If $0\ra M^\prime\ra M\ra M^{\prime\prime}\ra 0$
is an exact sequence in $\s$, then $0\ra T\bom_R M^\prime\ra T
\bom_R M\ra T\bom_R M^{\prime\prime}\ra 0$ is an exact sequence in
$\es_{T}$. Thus, $f(M)=f(M^\prime)+f(M^{\prime\prime})$, and so
there exists a unique group homomorphism $S_0(\varphi):S_0(R)\ra
S_0(T)$ such that $S_0(\varphi)([M])=[T\bom_R M]$, for all $M\in
\s$.

In particular, the inclusion map $\lambda:R\ra R[x]$ induces the group homomorphism $S_0(\lambda):S_0(R)\ra S_0(R[x])$ by
$S_0(\lambda)([M])=[M[x]]$, for all $M\in \s$.

\vskip2mm Let $T$ be a ring and let $u\in T$ be neither a
zero-divisor nor a unit. Write $\ol{T}=T/(u)$. For $M\in \es_T$,
there is an exact sequence $0\ra A\ra P\ra M\ra 0$, where $P$ is
finitely generated free over $T$ and $A\in\es_T$. Set $\ol{N}=N/uN$
for a $T$-module $N$. Thus Proposition \ref{p202} gives
$\ol{A},\ol{P}\in \es_{\ol{T}}$. Define $g:\es_T\ra S_0(\ol{T})$ by
$g(M)=[\ol{P}]-[\ol{A}]$.

\begin{lem}\label{l503}
{\rm
Assume $T$ and $g$ are as above. Then $g$ determines a group homomorphism $g_{\pi}$ from $S_0(T)$ to $S_0(\ol{T})$.}
\end{lem}

\no{\bf Proof.}\ \ First, we claim that $g$ is a well-defined map.
Indeed, if $0\ra B\ra Q\ra M\ra 0$ is another exact sequence over
$T$, where $Q$ is finitely generated free and $B\in\es_T$, then by
Schanuel's Lemma, $0\ra A\ra P\bigoplus B\ra Q\ra 0$ is an exact
sequence in $\es_T$. Since $u$ is regular on $Q,B,P,A$, the sequence
$0\ra \ol{A}\ra \ol{P}\bigoplus \ol{B}\ra \ol{Q}\ra 0$ is exact in
$\es_{\ol{T}}$. Then $[\ol{P}\bigoplus \ol{B}]=[\ol{A}]+[\ol{Q}]$ in
$S_0(\ol{T})$, i.e., $[\ol{P}]-[\ol{A}]=[\ol{Q}]-[\ol{B}]$.
Moreover, if $0\ra M_1\ra M\ra M_2\ra 0$ is an exact sequence in
$\es_T$, then we can construct a commutative $3\times 3$ diagram
with exact columns and rows:
$$\xymatrix@R=16pt@C=30pt{
&0\ar[d]&0\ar[d]&0\ar[d]\\
0\ar[r]&A_1\ar[d]\ar[r]&A\ar[d]\ar[r]&A_2\ar[d]\ar[r]&0\\
0\ar[r]&P_1\ar[d]\ar[r]&P\ar[d]\ar[r]&P_2\ar[d]\ar[r]&0\\
0\ar[r]&M_1\ar[d]\ar[r]&M\ar[d]\ar[r]&M_2\ar[d]\ar[r]&0\\
&0&0&0}$$ where $P,P_1,P_2$ are finitely generated free $T$-modules
and $A,A_1,A_2\in\es_T$. Also, since $u$ is regular on $P,P_1,P_2$,
the sequence $0\ra \ol{P_1}\ra \ol{P}\ra \ol{P_2}\ra 0$ is exact in
$\es_{\ol{T}}$, and so $[\ol{P}]=[\ol{P_1}]+[\ol{P_2}]$ in
$S_0(\ol{T})$. Similarly, $[\ol{A}]=[\ol{A_1}]+[\ol{A_2}]$ in
$S_0(\ol{R})$. Thus
$$g(M_1)+g(M_2)=[\ol{P_1}]-[\ol{A_1}]+[\ol{P_2}]-[\ol{A_2}]= [\ol{P}]-[\ol{A}]=g(M).$$
Therefore, there exists a unique group homomorphism
$g_{\pi}:S_0(T)\ra S_0(\ol{T})$ with
$g_{\pi}([M])=g(M)=[\ol{P}]-[\ol{A}]$ for all
$M\in\es_T$.\hfill$\Box$

\begin{rem}\label{r504}
{\rm
In the proof of Lemma \ref{l503}, we see that if $u$ is regular on $M$, then $g_{\pi}([M])=[\ol{M}]$.}
\end{rem}

The next result plays an important role in the proof of our main theorem in this section.

\begin{thm}\label{t505}
{\rm
For any ring $R$, $S_0(\lambda):S_0(R)\ra S_0(R[x])$ is an isomorphism.}
\end{thm}

\no{\bf Proof.}\ \ If $N\in\es_{R[x]}$, then there is an exact
sequence of $R[x]$-modules $0\ra K\ra F[x]\ra N\ra 0$, where $F$ is
a finitely generated free $R$-module and $K\in\es_{R[x]}$. Thus
$g_{\pi}([N])=[F]-[K/xK]$, where $g_{\pi}:S_0(R[x])\ra S_0(R)$ is
the homomorphism given in Lemma \ref{l503} for the case $T=R[x]$ and
$u=x$.

It is routine to verify that $g_{\pi}S_0(\lambda)={\bf 1}_{S_0(R)}$,
and so $S_0(\lambda)$ is a monomorphism. To complete the proof, we
need only show that $S_0(\lambda)$ is an epimorphism. For
$N\in\es_{R[x]}$, there are exact sequences of $R[x]$-modules $0\ra
N_0\ra P[x]\ra N\ra 0$ and $0\ra N_1\ra F[x]\ra N_0\ra 0$, where $P$
and $F$ are finitely generated free $R$-modules and
$N_0,N_1\in\es_{R[x]}$. Thus, Lemma \ref{l312} says that there
exists an exact sequence $0\ra A[x]\ra B[x]\ra N_1\ra 0$ in
$\es_{R[x]}$, where $A,B\in\s$. Therefore in $S_0(R[x])$, we have
$$[N]=[P[x]]-[N_0]=[P[x]]-[F[x]]+[N_1]=[P[x]]-[F[x]]+[B[x]]-[A[x]].$$
Therefore, $[N]=S_0(\lambda)([P]-[F]+[B]-[A])$, which implies that $S_0(\lambda)$ is epimorphic.\hfill$\Box$

\vskip2mm
Our next result shows that the property (C) is also preserved under polynomial extensions.

\begin{pro}\label{p506}
{\rm
If $R$ has the property (C), then so does $R[x]$.}
\end{pro}

\no{\bf Proof.}\ \ The proof of this proposition is similar to that given in Theorem \ref{t313}.\hfill$\Box$

From the above results we get the following main theorem in this section.

\begin{thm}\label{t507}
{\rm
If $R$ has the property (C), then $R$ is $K_0$-regular.}
\end{thm}

\no{\bf Proof.}\ \ It suffices, by induction, to show that the functorial map $K_0(\lambda):K_0(R)\ra K_0(R[x])$ is an isomorphism.
But this follows immediately from Proposition \ref{p506}, Proposition \ref{p501}, Theorem \ref{t505}, and the following
commutative diagram:
$$\xymatrix@R=20pt@C=40pt{
K_0(R)\ar[d]_{\theta_R}\ar[r]^{K_0(\lambda)}&K_0(R[x])\ar[d]^{\theta_{R[x]}}\\
S_0(R)\ar[r]^{S_0(\lambda)}&S_0(R[x]).}$$
\hfill$\Box$

\begin{cor}\label{c508}
{\rm
If $\spd(R)<\infty$, then $R$ is $K_0$-regular.}
\end{cor}

\begin{cor}\label{c509}
{\rm
If $\wgl(R)<\infty$, then $R$ is $K_0$-regular.}
\end{cor}

\vskip4mm
\bc{\bf Acknowledgements}\ec

\vskip3mm The authors would like to thank referee for several
valuable suggestions. This work was supported by NSFC (No. 11171240)
and the Specialized Research Fund for the Doctoral Program of Higher
Education (No. 20125134110002).

\end{document}